\documentclass{tran-l}

\usepackage[mathscr]{eucal}
\usepackage{amscd}
\usepackage{mathrsfs}
\usepackage{amsfonts}
\usepackage{amsmath, amsthm, amssymb}
\usepackage{latexsym}
\usepackage{calrsfs}
\usepackage[dvips]{graphics}
\usepackage{xcolor}

\usepackage{lipsum} 
\usepackage{cite}
\usepackage[sc]{mathpazo} 
\usepackage[T1]{fontenc} 
\usepackage{microtype} 

\usepackage[hmarginratio=1:1,top=30mm]{geometry} 
\usepackage{hyperref} 

\newtheorem{thm}{Theorem}[section]

\newtheorem{lem}[thm]{Lemma}
\newtheorem{fact}[thm]{Fact}

\newtheorem{exa}[thm]{Example}
\numberwithin{equation}{section}

\begin{document}

\title[Commutativity criteria for prime rings]{ Commutativity criteria for prime rings with involution via pairs of endomorphisms}

\author{G. S. Sandhu$^1$, G. Gudwani$^{2}$, M. El Hamdoui$^{3}$}

\address{$^1$ Department of Mathematics, 
Patel Memorial National College, Rajpura 140401, India}
\email{gurninder\_rs@pbi.ac.in}
\address{$^2$ School of Basic, Applied and Bio Sciences, RIMT University, Mandi Gobindgarh 147301, Punjab, India}
\email{pmn.geetika@gmail.com}
\address{$^3$ Faculty of Sciences and Technologies (FST), Sidi Mohammed Ben Abdellah University, Fes 30000, Morocco}
\email{mathsup2011@gmail.com}

%




\begin{abstract}
The aim of this article is to investigate central-valued identities involving pairs of endomorphisms on prime rings equipped with an involution of the second kind. Extending the recent contributions of Mir et al. (2020) and Boua et al. (2024), we establish several new commutativity criteria for such rings in the presence of two distinct nontrivial endomorphisms. Our approach provides a unified technique that covers multiple classes of $\ast$-identities and yields generalizations of earlier single-endomorphism results. Moreover, explicit counterexamples are constructed to demonstrate the necessity of the hypotheses on primeness and on the nature of the involution.
\end{abstract}

\maketitle

\section{Introduction}

Throughout this paper, $R$ denotes an associative ring, and $Z(R)$ represents its center. For any $x, y \in R$, the commutator and anti-commutator are defined respectively by:
	$
	[x, y] = xy - yx \quad \text{and} \quad x \circ y = xy + yx.
	$
	A ring $R$ is called \emph{prime} if $aRb = \{0\}$ implies $ a = 0 $ or $ b = 0 $. For a positive integer $ n $, $ R $ is $ n $-torsion free if $ nx = 0 $ for $ x \in R $ implies $ x = 0 $.
	
	An \emph{involution} on $ R $ is an anti-automorphism $ \ast \colon R \to R $ satisfying $ (x^\ast)^\ast = x $ for all $ x \in R $. An element $ x \in R $ is \emph{symmetric} if $ x^\ast = x $, and is \emph{skew-symmetric} if $ x^\ast = -x $. The sets of symmetric and skew-symmetric elements are denoted by $ H(R) $ and $ S(R) $, respectively. An involution $ \ast $ is of the \emph{first kind} if $ Z(R) \subseteq H(R) $; otherwise, it is of the \emph{second kind}, in which case $ S(R) \cap Z(R) \neq \{0\} $.
	
	\medskip
	
	The study of commutativity in prime rings equipped with non-identity endomorphisms has evolved into a vibrant research area \cite{Deng1996}. Progress in this field hinges on analyzing diverse conditions and structural constraints that induce commutativity. Recent works have explored the interplay between prime ring structures, endomorphism properties, and algebraic behaviors of their elements; see, e.g., \cite{Aharssi23, Boua24, Idrissi22, Mamouni20, Mir20}.
	
	Historically, Divinsky \cite{Divinsky1958} pioneered this direction by investigating automorphisms. Specifically, he proved that no non-commutative simple Artinian ring admits a non-trivial automorphism $ f $ satisfying $ [f(x), x] = 0 $ for all $ x \in R $. Earlier, Posner \cite{posner1957} studied derivations in prime rings, showing that if a prime ring $ R $ has a nonzero derivation $ d $ with $ [d(x), x] \in Z(R) $ for all $ x \in R $, then $ R $ must be commutative. Building on the foundational work of Divinsky and Posner, the study of endomorphisms and derivations in ring theory has flourished, with significant contributions from algebraists worldwide. In \cite{Ashraf2001}, authors proved that a prime ring $R$ essentially commutative if it possesses a derivation $d$ satisfying the condition: $d(xy)-xy\in Z(R)$ or $d(xy)-yx\in Z(R)$ on a nonzero two-sided ideal of $R.$ In 2007, Ashraf et al. \cite{Ashraf2007} extended these results by taking generalized derivation instead of derivation. Specifically, they showed that the prime ring $R$ must be commutative if there exists a generalized derivation $F$ associated with a nonzero derivation $d$, which satisfies any of the identities: $F(xy)-xy\in Z(R),$ $F(xy)+xy\in Z(R),$ $F(xy)-yx\in Z(R),$ $F(xy)+yx\in Z(R),$ $F(x)F(y)-xy\in Z(R),$ $F(x)F(y)+xy\in Z(R)$ on a nonzero two-sided ideal of $R.$ Recently, Mamouni et al. \cite{Mamouni2018} considered more general situations by taking $\ast$-differential identities involving generalized derivations in prime rings with involution of the second kind. More precisely, the authors investigated that an involutive prime ring $(R,\ast)$ will be commutative if there exists a generalized derivation $F$ associated with a nonzero derivation $d$ satisfying any one of the following $\ast$-differential identity: $F(xx^{\ast})-xx^{\ast}\in Z(R),$ $F(xx^{\ast})+xx^{\ast}\in Z(R),$ $F(xx^{\ast})-x^{\ast}x\in Z(R),$ $F(xx^{\ast})+x^{\ast}x\in Z(R),$ $F(x)F(x^{\ast})-xx^{\ast}\in Z(R),$ $F(x)F(x^{\ast})+xx^{\ast}\in Z(R).$ Apart from this, there are more interesting studies available with these settings, for instance, see \cite{Ali2015}, \cite{Bharat2021}, \cite{Zem19}. 
	
	In the present work, however, we shift our focus to the investigation of endomorphisms of prime rings and their recent developments.


 
\section{Main Results}
Recently, Mir \textit{et al.} \cite{Mir20} explored the commutativity of prime rings under the influence of a nontrivial endomorphism $g$, considering the conditions:  
$g(xy)+xy \in Z(R),~g(xy)-xy \in Z(R),~ g(xy)+yx \in Z(R), g(xy) -yx \in Z(R)$  
for all $x, y \in R$. Additionally, the authors extended these results to rings equipped with an involution $\ast$ of the second kind, proving that a 2-torsion free prime ring $R$ admitting such an involution and a nontrivial endomorphism $g$ must be commutative if any of the following $\ast$-identities hold:  
$
g(xx^{\ast})+xx^{\ast} \in Z(R),~g(xx^{\ast})-xx^{\ast} \in Z(R),~ g(xx^{\ast})+x^{\ast}x \in Z(R),~ g(xx^{\ast}) - x^{\ast}x \in Z(R).
$  

In a further extension, Boua \textit{et al.} \cite{Boua24} examined prime rings endowed with two endomorphisms $g_1$ and $g_2$, establishing commutativity under the conditions:  
$g_1(x)g_2(y) + xy \in Z(R),~ g_1(x)g_2(y) - xy \in Z(R),~g_1(x)g_2(y) + yx \in Z(R),~ g_1(x)g_2(y) - yx \in Z(R).
$  

A natural question arises: if $R$ is a prime ring with an involution $\ast$ (of the second kind) and two endomorphisms $g_1, g_2$, can the commutativity of $R$ be deduced from any one of the following $\ast$-identities?
	\begin{itemize}
		\item[$(A_{1})$] $g_1(x)g_2(x^{\ast})-xx^{\ast} \in Z(R)$ for all $x\in R;$
		\item[$(A_{2})$] $g_1(x)g_2(x^{\ast}) + xx^{\ast} \in Z(R)$ for all $x\in R;$ 
		\item[$(A_{3})$] $g_1(x)g_2(x^{\ast})- x^{\ast}x \in Z(R)$ for all $x\in R;$
		\item[$(A_{4})$] $g_1(x)g_2(x^{\ast})+ x^{\ast}x \in Z(R)$ for all $x\in R.$
	\end{itemize}  

In this article, we address this problem by developing a unified approach applicable to all such scenarios. Our methodology not only resolves the posed question but also provides a streamlined framework for related investigations in the theory of prime rings with involution.
\par
The following are some preliminary facts that are crucial in the development of main results in this study. 

\begin{fact}{\cite[Theorem 3.4]{Aharssi23}}\label{Fact4}
	Let $R$ be a prime ring with involution $\ast$ of the second kind. Then the following conditions are equivalent:
	\begin{itemize}
		\item[(i)] $h^{2}\in Z(R)$ for all $h\in H(R);$
		\item[(ii)] $s^{2}\in Z(R)$ for all $s\in S(R);$
		\item[(iii)] $R$ is commutative. 
	\end{itemize}
\end{fact}

\begin{fact}{\cite[Fact 2]{Idrissi22}}\label{Fact3}
	Let $R$ be a 2-torsion free prime ring with involution $\ast$ of the second kind. If $[x\circ x^{\ast},a]\in Z(R)$ for all $x\in R,$ then $a\in Z(R)$.
\end{fact}
\begin{fact}{\cite[Lemma 2.1]{Nejjar17}}\label{Fact2}
	Let $R$ be a prime ring with involution $\ast$ of the second kind. Then $[x,x^{\ast}]\in Z(R)$ for all $x\in R$ if and only if $R$ is commutative.
\end{fact}
Throughout this paper, by a nontrivial endomorphism, we mean an endomorphism which is not identity on the given ring.

\begin{lem}\label{LEM-0}
	Let $R$ be a 2-torsion free prime ring with involution $\ast$ of the second kind. If $R$ admits nontrivial endomorphisms $g_{1},g_{2}$ satisfying any of $(A_{1})$-$(A_{4})$, then either $R$ is commutative or $g_{1}(Z(R))\subseteq Z(R)$ and $g_{2}(Z(R))\subseteq Z(R).$
	\begin{proof}
		Let us suppose that $(A_{1})$ holds, so we have
		\begin{equation}\label{A-1}
			g_{1}(x)g_{2}(x^{\ast})-xx^{\ast}\in Z(R),~\forall~x\in R.	
		\end{equation}
		Linearizing this relation, we get
		\begin{equation}\label{A-1.1}
			g_{1}(x)g_{2}(y^{\ast})+g_{1}(y)g_{2}(x^{\ast})-xy^{\ast}-yx^{\ast}\in Z(R),~\forall~x,y\in R.	
		\end{equation}
		For some $h\in H(R)\cap Z(R),$ taking $x=h$ in the relation (\ref{A-1}), we obtain $g_{1}(h)g_{2}(h)-h^{2}\in Z(R),$ which forces $g_{1}(h)g_{2}(h)\in Z(R).$
		Replacing $x$ by $xh$ in (\ref{A-1}), we have
		\[
		g_{1}(x)g_{1}(h)g_{2}(h)g_{2}(x^{\ast})- xx^{\ast}h^{2}\in Z(R),~\forall~x\in R,h\in H(R)\cap Z(R).	
		\]
		It implies that
		\begin{equation}\label{A-2}
			g_{1}(x)g_{2}(x^{\ast})g_{1}(h)g_{2}(h)- xx^{\ast}h^{2}\in Z(R),~\forall~x\in R,h\in H(R)\cap Z(R).	
		\end{equation}
		Multiplying (\ref{A-1}) by $h^{2}$ and then subtracting from (\ref{A-2}), we get
		\begin{equation*}
			g_{1}(x)g_{2}(x^{\ast})(g_{1}(h)g_{2}(h)-h^{2})\in Z(R),~\forall~x\in R,h\in H(R)\cap Z(R).	
		\end{equation*}
		Using primeness of $R$ we obtain that either $g_{1}(x)g_{2}(x^{\ast})\in Z(R)$ for all $x\in R$ or $g_{1}(h)g_{2}(h)=h^{2}$ for all $h\in H(R)\cap Z(R).$ In the first case it follows from our hypothesis that $xx^{\ast}\in Z(R)$ for all $x\in R.$ Replacing $x$ by $x^{\ast}$ in the last relation, we find $x^{\ast}x\in Z(R)$ for all $x\in R.$ It forces that $[x,x^{\ast}]\in Z(R)$ for all $x\in R,$ hence $R$ is commutative by Fact \ref{Fact2}.
		\par
		Now, let us consider the latter case, i.e., $g_{1}(h)g_{2}(h)=h^{2}$ for all $h\in H(R)\cap Z(R).$
		Substituting $y=h,$ where $0\neq h\in H(R)\cap Z(R),$ in (\ref{A-1.1}) equation and we get
		\[
		g_{1}(x)g_{2}(h)+g_{1}(h)g_{2}(x^{\ast})-(x+x^{\ast})h\in Z(R),~\forall~x\in R.	
		\]
		Replacing $x$ by $xh$ in the last obtained expression, we may infer that
		\[
		g_{1}(x)g_{1}(h)g_{2}(h)+g_{1}(h)g_{2}(h)g_{2}(x^{\ast})
		-(x+x^{\ast})h^{2}\in Z(R),~\forall~x\in R,~0\neq h\in H(R)\cap Z(R).	
		\]
		Thus, it gives
		\[
		\bigg(g_{1}(x)+g_{2}(x^{\ast})
		-(x+x^{\ast})\bigg)h^{2}\in Z(R),~\forall~x\in R,~0\neq h\in H(R)\cap Z(R).	
		\]
		Since $0\neq h,$ so we conclude that
		\begin{equation}\label{A-4}
			g_{1}(x)+g_{2}(x^{\ast})-(x+x^{\ast})\in Z(R),~\forall~x\in R.
		\end{equation}
		Replacing $x$ by $x-x^{\ast}$ in (\ref{A-4}), we get
			\begin{equation}\label{A-5}
				g_{1}(x-x^{\ast})+g_{2}(x-x^{\ast})\in Z(R),~\forall~x\in R.
			\end{equation}
			For $c\in Z(R),$ we see from (\ref{A-4}) that
		\begin{equation}\label{A-6}
			[g_{1}(x),g_{2}(c)]=[x+x^{\ast},g_{2}(c)],~\forall~x\in R.
		\end{equation}
		Replacing $x$ by $x^{\ast}$ in (\ref{A-6}), we have
		\begin{equation}\label{A-7}
			[g_{1}(x^{\ast}),g_{2}(c)]=[x+x^{\ast},g_{2}(c)],~\forall~x\in R.
		\end{equation}
		Subtracting (\ref{A-6}) from (\ref{A-7}), we conclude that
		\begin{equation}\label{A-8}
			[g_{1}(x-x^{\ast}),g_{2}(c)]=0,~\forall~x\in R.
		\end{equation}
		Substituting $x-x^{\ast}$ instead of $x$ in (\ref{A-1}), we can see that
		\[
		[g_{1}(x-x^{\ast})g_{2}(x-x^{\ast})+(x-x^{\ast})^{2},g_{2}(c)]=0,~\forall~x\in R.
		\]
		In view of (\ref{A-8}), it follows that 
		\begin{equation}\label{A-9}
			[x^{2}-xx^{\ast}-x^{\ast}x+(x^{\ast})^{2},g_{2}(c)]=0,~\forall~x\in R.
		\end{equation}
		Replacing $x$ by $xk_{c}$ in the above equation, where $0\neq k_{c}\in S(R)\cap Z(R),$ we have
		\[
		[x^{2}+xx^{\ast}+x^{\ast}x+(x^{\ast})^{2},g_{2}(c)]k_{c}^{2}=0,~\forall~x\in R.
		\]
		As $k_{c}\neq 0,$ it forces
		\begin{equation}\label{A-10}
			[x^{2}+xx^{\ast}+x^{\ast}x+(x^{\ast})^{2},g_{2}(c)]=0,~\forall~x\in R.
		\end{equation}
		Subtracting (\ref{A-9}) from (\ref{A-10}), we may infer that $[x\circ x^{\ast},g_{2}(c)]=0$ for all $x\in R.$ Fact \ref{Fact3} yields that $g_{2}(Z(R))\subseteq Z(R).$ Now from (\ref{A-4}), it is straightforward that $g_{1}(Z(R))\subseteq Z(R).$
\\
With similar argument, we can observe the conclusion for other situations $(A_{2})-(A_{4}).$ 
	\end{proof} 
\end{lem}

\begin{lem}\label{LEM-1}
Let $R$ be a 2-torsion free prime ring with involution $\ast$ of the second kind. If $R$ admits nontrivial endomorphisms $g_{1},g_{2}$ satisfying $(A_{1})$ or $(A_{3})$ with $g_{1}(h)=h$ for all $h\in H(R)\cap Z(R),$ then $R$ is commutative.
	\begin{proof} Let assume that $(A_{1})$ holds.
		Suppose $R$ is not commutative. By hypothesis, we have
		\begin{equation}\label{A-0}
			g_{1}(x)g_{2}(x^{\ast})- xx^{\ast}\in Z(R),~\forall~x\in R.	
		\end{equation}
		Let us first polarize this relation to get
		\begin{equation}\label{A-3}
			g_{1}(x)g_{2}(y^{\ast})+g_{1}(y)g_{2}(x^{\ast})-xy^{\ast}- yx^{\ast}\in Z(R),~\forall~x,y\in R.	
		\end{equation}
		Substituting $xh$ for $x$ in (\ref{A-0}), we have
		\begin{equation}\label{A-1a}
			g_{1}(x)g_{2}(x^{\ast})hg_{2}(h)- xx^{\ast}h^{2}\in Z(R),~\forall~x\in R.	
		\end{equation}
		Moreover, for $h\in H(R)\cap Z(R),$ taking $x=h$ in (\ref{A-0}), we obtain $g_{1}(h)g_{2}(h)-h^{2}\in Z(R).$ Thus, we also get $hg_{2}(h)\in Z(R).$
		\\
		Multiplying (\ref{A-0}) by $hg_{2}(h),$ we find
		\begin{equation}\label{A-1b}
			g_{1}(x)g_{2}(x^{\ast})hg_{2}(h)- xx^{\ast}hg_{2}(h)\in Z(R),~\forall~x\in R.	
		\end{equation}
		 Comparing (\ref{A-1a}) and (\ref{A-1b}), it follows that
		 \[
		 xx^{\ast}h(g_{2}(h)-h)\in Z(R),~\forall~x\in R.
		 \]
		 It forces that either $xx^{\ast}\in Z(R)$ for all $x\in R$ or $h(g_{2}(h)-h)=0.$ In the first case we get $R$ commutative, and we are done. In the latter case, as $h\neq 0,$ we find $g_{2}(h)=h$ for all $h\in H(R)\cap Z(R).$
		 \par
		This means that $g_{1}(h)=g_{2}(h)=h$ for all $h\in H(R)\cap Z(R).$ Replacing $h$ by $k_{c}^{2},$ where $k_{c}\in S(R)\cap Z(R),$ we get $g_{1}(k_{c})^{2}=k_{c}^{2}.$ It takes us to  $(g_{1}(k_{c})-k_{c})(g_{1}(k_{c})+k_{c})=0.$ Invoking Lemma \ref{LEM-0}, we obtain that either $g_{1}(k_{c})=k_{c}$ or $g_{1}(k_{c})=-k_{c}.$ 
		\par Similarly for all $k_{c}\in S(R)\cap Z(R),$ we get $g_{2}(k_{c})=k_{c}$ or $g_{2}(k_{c})=-k_{c}.$ Thus, the following cases arise:
		\begin{itemize}
			\item[(i)] $g_{1}(k_{c})=k_{c}$ and $g_{2}(k_{c})=k_{c}$ for all $k_{c}\in S(R)\cap Z(R);$
			\item[(ii)] $g_{1}(k_{c})=-k_{c}$ and $g_{2}(k_{c})=-k_{c}$ for all $k_{c}\in S(R)\cap Z(R);$
			\item[(iii)] $g_{1}(k_{c})=-k_{c}$ and $g_{2}(k_{c})=k_{c}$ for all $k_{c}\in S(R)\cap Z(R);$
			\item[(iv)] $g_{1}(k_{c})=k_{c}$ and $g_{2}(k_{c})=-k_{c}$ for all $k_{c}\in S(R)\cap Z(R).$
		\end{itemize}
		Let us discuss these cases in detail and see their respective outcomes.
		\\
		(i) Suppose that $g_{1}(k_{c})=k_{c}$ and $g_{2}(k_{c})=k_{c}$ for all $k_{c}\in S(R)\cap Z(R).$ Replacing $y$ by $yk_{c}$ in (\ref{A-3}), where $0\neq k_{c}\in S(R)\cap Z(R),$ we obtain
		\begin{equation}\label{A-3a}
			\bigg(-g_{1}(x)g_{2}(y^{\ast})+g_{1}(y)g_{2}(x^{\ast})+xy^{\ast}- yx^{\ast}\bigg)k_{c}\in Z(R),~\forall~x,y\in R.	
		\end{equation}
		Multiplying (\ref{A-3}) by $k_{c},$ we get
		\begin{equation}\label{A-3b}
			\bigg(g_{1}(x)g_{2}(y^{\ast})+g_{1}(y)g_{2}(x^{\ast})-xy^{\ast}- yx^{\ast}\bigg)k_{c}\in Z(R),~\forall~x,y\in R.	
		\end{equation}
		Adding (\ref{A-3a}) and (\ref{A-3b}), we get $g_{1}(y)g_{2}(x^{\ast})-yx^{\ast}\in Z(R)$ for all $x,y\in R.$ Putting $x=x^{\ast},$ we get $g_{1}(y)g_{2}(x)-yx\in Z(R)$ for all $x,y\in R.$ Thus the conclusion follows from Corollary 2 of \cite{Boua24}.
		\vspace{0.5cm}
		
		(ii) Assume that $g_{1}(k_{c})=-k_{c}$ and $g_{2}(k_{c})=-k_{c}$ for all $k_{c}\in S(R)\cap Z(R).$ Substituting $yk_{c}$ for $y$ in (\ref{A-3}), where $0\neq k_{c}\in S(R)\cap Z(R)$, we have 
		\[
		\Big(g_{1}(x)g_{2}(y^{\ast})-g_{1}(y)g_{2}(x^{\ast})+xy^{\ast}- yx^{\ast}\Big)k_{c}\in Z(R),~\forall~x,y\in R.	
		\]
		Since $k_{c}$ is nonzero, by primeness of $R,$ we find that
		\begin{equation}\label{AC-2}
			g_{1}(x)g_{2}(y^{\ast})-g_{1}(y)g_{2}(x^{\ast})+xy^{\ast}- yx^{\ast}\in Z(R),~\forall~x,y\in R.
		\end{equation}
		Combining (\ref{AA-3}) and (\ref{AB-2}) to conclude that
		$g_{1}(x)g_{2}(y^{\ast})-yx^{\ast}\in Z(R)$ for all $x,y\in R.$
		Taking $x=h,$ where $0\neq h\in H(R)\cap Z(R),$ we get $g_{2}(y^{\ast})-y\in Z(R)$ for all $y\in R.$ If we replace $y$ by 
		$y^{\ast}$ in this relation, then we find an identical expression to Eq. (15) of \cite{Mir20}, repeating the same reasoning, we get our conclusion.
		\vspace{0.5cm}
		
		(iii) Suppose that $g_{1}(k_{c})=-k_{c}$ and $g_{2}(k_{c})=k_{c}$ for all $k_{c}\in S(R)\cap Z(R).$
		Replacing $y$ by $yk_{c}$ in (\ref{A-3}), where $0\neq k_{c}\in S(R)\cap Z(R),$ we obtain
		\begin{equation}\label{A-3c}
			\bigg(g_{1}(x)g_{2}(y^{\ast})+g_{1}(y)g_{2}(x^{\ast})+xy^{\ast}- yx^{\ast}\bigg)k_{c}\in Z(R),~\forall~x,y\in R.	
		\end{equation}
		Multiplying (\ref{A-3}) by $k_{c}$ and adding with (\ref{A-3c}), it yields that $yx^{\ast}\in Z(R)$ for all $x,y\in R.$ It easily implies the commutativity of $R.$
		\vspace{0.5cm}
		
		(iv) Finally, we consider that $g_{1}(k_{c})=k_{c}$ and $g_{2}(k_{c})=-k_{c}$ for all $k_{c}\in S(R)\cap Z(R).$ It can be seen that this case follows in the same manner as (iii), so we omit its proof. 
		\\ It completes the proof.
		\par It can be verified that the same conclusion follows from condition $(A_{3})$ by an argument parallel to the above, with only minor modifications. Hence, the proof is omitted.
	\end{proof}
\end{lem}
\begin{thm}\label{thm-3}
	Let $R$ be a 2-torsion free prime ring with involution $\ast$ of the second kind. If $R$ admits nontrivial endomorphisms $g_{1},g_{2}$ satisfying $(A_{1})$ or $(A_{3}),$
	then $R$ is a commutative.
\begin{proof}
	Suppose that $R$ is not commutative. Then by Lemma \ref{LEM-0}, we know that $g_{1}$ and $g_{2}$ both are center-preserving. We shall now use this fact implicitly in order to proceed further and also recall some relations obtain in Lemma \ref{LEM-0}.
Replacing $x$ by $xh$ in (\ref{A-4}), we have
		\begin{equation}\label{A-11}
			g_{1}(x)g_{1}(h)+g_{2}(x^{\ast})g_{2}(h)-(x+x^{\ast})h\in Z(R),~\forall~x\in R.
		\end{equation}
		Multiplying (\ref{A-4}) by $h$ and subtracting from (\ref{A-11}) in order to find
		\begin{equation}\label{A-12}
			g_{1}(x)(g_{1}(h)-h)+g_{2}(x^{\ast})(g_{2}(h)-h)\in Z(R),~\forall~x\in R.
		\end{equation}
		Replacing $x$ by $xk_{c}$ in the above equation, we get
		\begin{equation}\label{A-13}
			g_{1}(k_{c})g_{1}(x)(g_{1}(h)-h)-g_{2}(k_{c})g_{2}(x^{\ast})(g_{2}(h)-h)\in Z(R),~\forall~x\in R,~k_{c}\in S(R)\cap Z(R).
		\end{equation}
		Multiplying (\ref{A-12}) by $g_{2}(k_{c})$ and compare with (\ref{A-13}), we obtain that
		\[
		(g_{2}(k_{c})+g_{1}(k_{c}))g_{1}(x)(g_{1}(h)-h)\in Z(R),~\forall~x\in R,~k_{c}\in S(R)\cap Z(R).
		\]
		In view of primeness of $R$, we have the following cases:
		\begin{itemize}
			\item[(a)] $g_{1}(x)\in Z(R)$ for all $x\in R;$
			\item[(b)] $g_{2}(k_{c})+g_{1}(k_{c})=0$ for all $k_{c}\in S(R)\cap Z(R);$ 
			\item[(c)] $g_{1}(h)=h$ for all $h\in H(R)\cap Z(R).$
		\end{itemize}
		\noindent \textbf{Case (a).} If $g_{1}(x) \in Z(R)$, then multiplying (\ref{A-5}) by $g_{1}((x-x^{\ast})^{\ast})$ in order to see that
		\[
		g_{1}((x-x^{\ast})^{\ast})g_{1}(x-x^{\ast})+g_{1}((x-x^{\ast})^{\ast})g_{2}(x-x^{\ast})\in Z(R),~\forall~x\in R.
		\]
		Using our initial hypothesis, it yields that
		\[
		(x-x^{\ast})^{\ast}(x-x^{\ast})\in Z(R),~\forall~x\in R.
		\]
		It means, we have $(x-x^{\ast})^{2}\in Z(R),~\forall~x\in R.$
		In particular, putting $x=s,$ where $s\in S(R),$ we see that $s^{2}\in Z(R),$ as $R$ is 2-torsion free.
		Now by Fact \ref{Fact4}, it forces $R$ to be commutative, which is a contradiction.\\
		\noindent \textbf{Case (b).} 
		Let us consider that $g_{2}(k_{c})=-g_{1}(k_{c})$ for all $k_{c}\in S(R)\cap Z(R).$ Replacing $x$ by $xk_{c}$ in (\ref{A-5}), we get
		\[
		(g_{1}(x+x^{\ast})-g_{2}(x+x^{\ast}))g_{1}(k_{c})\in Z(R),~\forall~x\in R.
		\]
		It forces that either $g_{1}(k_{c})=0$ or $g_{1}(x+x^{\ast})-g_{2}(x+x^{\ast})\in Z(R)$ for all $x\in R.$ Let us consider the latter case first. So we have
		\begin{equation}\label{A-14}
			g_{1}(x+x^{\ast})-g_{2}(x+x^{\ast})\in Z(R),~\forall~x\in R.
		\end{equation}
		Combining (\ref{A-5}) and (\ref{A-14}), we get
		\begin{equation}\label{A-15}
			g_{1}(x)-g_{2}(x^{\ast})\in Z(R),~\forall~x\in R.
		\end{equation}
		Now let us adding (\ref{A-4}) and (\ref{A-15}) in order to obtain
		\[
		2g_{1}(x)-(x+x^{\ast})\in Z(R),~\forall~x\in R.
		\]
		Changing $x$ into $x-x^{\ast}$ in the last relation, we get 
		\begin{equation}\label{A-16}
			g_{1}(x-x^{\ast})\in Z(R),~\forall~x\in R.
		\end{equation}
		Replacing $x$ by $xk_{c}$ in the above relation, we have
		$g_{1}(x+x^{\ast})g_{1}(k_{c})\in Z(R)$ for all $x\in R.$ By primeness of $R$ we have either 
		$g_{1}(x+x^{\ast})\in Z(R)$ for all $x\in R$ or $g_{1}(k_{c})=0$ for all $k_{c}\in S(R)\cap Z(R).$ Let us first assume that
		\begin{equation}\label{A-17}
			g_{1}(x+x^{\ast})\in Z(R),~\forall~x\in R.
		\end{equation}
		Jointly considering (\ref{A-16}) and (\ref{A-17}), we get that $g_{1}(x)\in Z(R)$ for all $x\in R,$ which is our case (a), so we are done. 
		
		On the other hand, if $g_{1}(k_{c})=0$ for all $k_{c}\in S(R)\cap Z(R)$; by substituting $xk_{c}$ for $x$ in (\ref{A-4}), where $0\neq k_{c}\in S(R)\cap Z(R),$ we can see that $x-x^{\ast}\in Z(R)$ for all $x\in R.$ It immediately takes us to that $R$ is commutative, a contradiction.\\ 
		\noindent \textbf{Case (c).} This case is done by Lemma \ref{LEM-1}. 
		\\Therefore, we conclude that $R$ is commutative, as desired.
		\par It can be verified that the same conclusion follows from condition $(A_{3})$ by an argument parallel to the above, with only minor modifications. Hence, the proof is omitted.
	\end{proof}
\end{thm}

\begin{lem}\label{LEM-2}
	Let $R$ be a 2-torsion free prime ring with involution $\ast$ of the second kind. If $R$ admits nontrivial endomorphisms $g_{1},g_{2}$ satisfying $(A_{2})$ or $(A_{4})$
	with $g_{1}(h)=h$ for all $h\in H(R)\cap Z(R),$ then $R$ is commutative.
	\begin{proof}
The proof proceeds almost in the same manner as that of Lemma \ref{LEM-1}, except for some variations in the mapping $g_{2}$, which we shall now examine.
Following in the same way we get that $g_{2}(h)=-h$ for all $h\in H(R)\cap Z(R)$. Moreover, for each $k_{c}\in S(R)\cap Z(R),$ either $g_{1}(k_{c})=k_{c}$ or $g_{1}(k_{c})=-k_{c}.$ But we can not explicitly have the same behavior of $g_{2}$ on central skew-symmetric elements. In order to obtain that, we need more reasoning. 
Polarizing the given hypothesis, we have
\begin{equation}\label{AA-3}
	g_{1}(x)g_{2}(y^{\ast})+g_{1}(y)g_{2}(x^{\ast})+xy^{\ast}+ yx^{\ast}\in Z(R),~\forall~x,y\in R.	
\end{equation}		
Taking $y=h$ in particular for $0\neq h\in H(R)\cap Z(R)$, we find that		
\begin{equation}\label{AA-3a}
	-g_{1}(x)+g_{2}(x^{\ast})+x+x^{\ast}\in Z(R),~\forall~x\in R.	
\end{equation}
First, let us suppose that $g_{1}(k_{c})=-k_{c}$ for all $k_{c}\in S(R)\cap Z(R).$ Replacing $x$ by $xk_{c}$ in (\ref{AA-3a}), one can verify that
\begin{equation}\label{AA-3b}
	g_{1}(x)k_{c}-g_{2}(x^{\ast})g_{2}(k_{c})+xk_{c}-x^{\ast}k_{c}\in Z(R),~\forall~x\in R.	
\end{equation}
Multiplying (\ref{AA-3a}) by $k_{c}$ and adding in (\ref{AA-3b}), we have $g_{2}(x^{\ast})(k_{c}-g_{2}(k_{c}))+2xk_{c}\in Z(R)$ for all $x\in R.$ It yields that
\begin{equation}\label{AA-3c}
	g_{2}(x)(k_{c}-g_{2}(k_{c}))+2x^{\ast}k_{c}\in Z(R),~\forall~x\in R.	
\end{equation}
Again substituting $xk_{c}$ for $x$ in (\ref{AA-3c}), we may infer that
\begin{equation}\label{AA-3d}
	g_{2}(x)g_{2}(k_{c})(k_{c}-g_{2}(k_{c}))-2x^{\ast}k_{c}^{2}\in Z(R),~\forall~x\in R.	
\end{equation}
Multiplying (\ref{AA-3c}) by $k_{c}$ and merging with (\ref{AA-3d}), we see that $g_{2}(x)(k_{c}+g_{2}(k_{c}))(k_{c}-g_{2}(k_{c}))\in Z(R)$ for all $x\in R$ and $k_{c}\in S(R)\cap Z(R).$ Since $g_{2}(Z(R))\subseteq Z(R)$  (see Lemma \ref{LEM-0}), by using primeness of $R$ we conclude that either $g_{2}(R)\subseteq Z(R)$ or $(k_{c}+g_{2}(k_{c}))(k_{c}-g_{2}(k_{c}))=0$ for all $k_{c}\in S(R)\cap Z(R).$ 
\par
Next, if $g_{1}(k_{c})=k_{c}$ for all $k_{c}\in S(R)\cap Z(R).$
Changing $x$ by $xk_{c}$ in (\ref{AA-3a}), one can verify that
\begin{equation}\label{AA-3e}
	-g_{1}(x)k_{c}-g_{2}(x^{\ast})g_{2}(k_{c})+xk_{c}-x^{\ast}k_{c}\in Z(R),~\forall~x\in R.	
\end{equation}
Multiplying (\ref{AA-3a}) by $k_{c}$ and adding with (\ref{AA-3e}), we examine that
\begin{equation}\label{AA-3f}
	-2g_{1}(x)k_{c}+g_{2}(x^{\ast})(k_{c}-g_{2}(k_{c}))+2xk_{c}\in Z(R),~\forall~x\in R.	
\end{equation}
Writing $xk_{c}$ for $x$ in (\ref{AA-3f}) and simplify it, we arrive at $g_{2}(x^{\ast})(g_{2}(k_{c})-k_{c})(g_{2}(k_{c})+k_{c})\in Z(R).$ Thence either $g_{2}(R)\subseteq Z(R)$ or $(g_{2}(k_{c})-k_{c})(g_{2}(k_{c})+k_{c})=0$ for all $k_{c}\in S(R)\cap Z(R).$

It shows that in each case we have the same conclusion, i.e., either $g_{2}(R)\subseteq Z(R)$ or $(g_{2}(k_{c})-k_{c})(g_{2}(k_{c})+k_{c})=0$ for all $k_{c}\in S(R)\cap Z(R).$ 

First let us discuss the case $g_{2}(R)\subseteq Z(R).$ From (\ref{AA-3a}), we observe that $-g_{1}(x)+x+x^{\ast}\in Z(R)$ for all $x\in R.$ For some $0\neq k_{c},$ putting $x=xk_{c}$ in it, we get	$-g_{1}(x)+x-x^{\ast}\in Z(R)$ for all $x\in R.$ Comparing last two expressions obtained, we find that $2x^{\ast}\in Z(R)$ for all $x\in R,$ it forces $R$ to be commutative, and hence we are done. On the other hand, we have the following cases:
\begin{itemize}
	\item[(i)] $g_{1}(k_{c})=k_{c}$ and $g_{2}(k_{c})=k_{c}$ for all $k_{c}\in S(R)\cap Z(R);$
	\item[(ii)] $g_{1}(k_{c})=-k_{c}$ and $g_{2}(k_{c})=-k_{c}$ for all $k_{c}\in S(R)\cap Z(R);$
	\item[(iii)] $g_{1}(k_{c})=-k_{c}$ and $g_{2}(k_{c})=k_{c}$ for all $k_{c}\in S(R)\cap Z(R);$
	\item[(iv)] $g_{1}(k_{c})=k_{c}$ and $g_{2}(k_{c})=-k_{c}$ for all $k_{c}\in S(R)\cap Z(R).$
\end{itemize}
(i) Let us suppose that $g_{1}(k_{c})=k_{c}$ and $g_{2}(k_{c})=k_{c}$ for all $k_{c}\in S(R)\cap Z(R).$ Placing $yk_{c}$ instead of $y$ in (\ref{AA-3}), where $0\neq k_{c}\in S(R)\cap Z(R)$, we have 
\[
	\Big(-g_{1}(x)g_{2}(y^{\ast})+g_{1}(y)g_{2}(x^{\ast})-xy^{\ast}+ yx^{\ast}\Big)k_{c}\in Z(R),~\forall~x,y\in R.	
\]
Since $k_{c}$ is nonzero, by primeness of $R,$ we find that
\begin{equation}\label{AB-1}
	-g_{1}(x)g_{2}(y^{\ast})+g_{1}(y)g_{2}(x^{\ast})-xy^{\ast}+ yx^{\ast}\in Z(R),~\forall~x,y\in R.
\end{equation}
Jointly considering (\ref{AA-3}) and (\ref{AB-1}) to conclude that
$g_{1}(y)g_{2}(x^{\ast})+yx^{\ast}\in Z(R)$ for all $x,y\in R.$ Replacing $x$ by $x^{\ast},$ we obtain $g_{1}(y)g_{2}(x)+yx\in Z(R)$ for all $x,y\in R.$ Hence, the conclusion is obtained by \cite[Theorem 2.4]{Boua24}.
\vspace{0.5cm}

(ii) Assume that $g_{1}(k_{c})=-k_{c}$ and $g_{2}(k_{c})=-k_{c}$ for all $k_{c}\in S(R)\cap Z(R).$ Substituting $yk_{c}$ for $y$ in (\ref{AA-3}), where $0\neq k_{c}\in S(R)\cap Z(R)$, we have 
\[
\Big(g_{1}(x)g_{2}(y^{\ast})-g_{1}(y)g_{2}(x^{\ast})-xy^{\ast}+ yx^{\ast}\Big)k_{c}\in Z(R),~\forall~x,y\in R.	
\]
Since $k_{c}$ is nonzero, by primeness of $R,$ we find that
\begin{equation}\label{AB-2}
	g_{1}(x)g_{2}(y^{\ast})-g_{1}(y)g_{2}(x^{\ast})-xy^{\ast}+ yx^{\ast}\in Z(R),~\forall~x,y\in R.
\end{equation}
Combining (\ref{AA-3}) and (\ref{AB-2}) to conclude that
$g_{1}(x)g_{2}(y^{\ast})+yx^{\ast}\in Z(R)$ for all $x,y\in R.$
Taking $x=h,$ where $0\neq h\in H(R)\cap Z(R),$ we get $g_{2}(y^{\ast})+y\in Z(R)$ for all $y\in R.$ If we replace $y$ by 
$y^{\ast}$ in this relation, then we find an identical expression to Eq. (19) of \cite{Mir20}, following in the same way, we get our conclusion.
\vspace{0.5cm}

(iii) Let us consider $g_{1}(k_{c})=-k_{c}$ and $g_{2}(k_{c})=k_{c}$ for all $k_{c}\in S(R)\cap Z(R).$ Substituting $yk_{c}$ for $y$ in (\ref{AA-3}), where $0\neq k_{c}\in S(R)\cap Z(R)$, we find that
\begin{equation}\label{AB-3}
	-g_{1}(x)g_{2}(y^{\ast})-g_{1}(y)g_{2}(x^{\ast})-xy^{\ast}+ yx^{\ast}\in Z(R),~\forall~x,y\in R.
\end{equation}
On adding this in (\ref{AA-3}), we get $yx^{\ast}\in Z(R)$ for all $x,y\in R.$ It forces commutativity of $R,$ which is a contradiction. Hence we are done.
\vspace{0.5cm}

(iv) finally, let us assume that $g_{1}(k_{c})=k_{c}$ and $g_{2}(k_{c})=-k_{c}$ for all $k_{c}\in S(R)\cap Z(R).$ This case can be easily seen in the same manner of Case (iii). 
\\ Therefore, $R$ is commutative.
\par It can be verified that the same conclusion follows from condition $(A_{4})$ by an argument parallel to the above, with only minor modifications. Hence, the proof is omitted.
\end{proof}
\end{lem}
\begin{thm}\label{thm-4}
	Let $R$ be a 2-torsion free prime ring with involution $\ast$ of the second kind. If $R$ admits nontrivial endomorphisms $g_{1},g_{2}$ satisfying $(A_{2}),$
	then $R$ is a commutative.
	\begin{proof}
		The proof follows on the same lines of Theorem \ref{thm-3}.
	\end{proof}
\end{thm}

The following example exhibits that the condition on involution $\ast$ is of the second kind in Theorem \ref{thm-3} and \ref{thm-4} is not superfluous.
\begin{exa}\label{example5.1}
	Let $R=M_{2}(\mathbb{Q})=\left\{ \left(\begin{array}{cc}
		a & b \\
		c & d
	\end{array} \right) |~a,b,c,d \in \mathbb{Q} \right\}$, where $\mathbb{Q}$ denotes the ring of rational numbers. Define the mappings $\ast, g:R\to R$ by
	$\left( \begin{array}{cc}
		a & b \\
		c & d
	\end{array} \right)^\ast= \left( \begin{array}{cc}
		d & -b \\
		-c & a
	\end{array} \right),$
	$g\left( \begin{array}{cc}
		a & b \\
		c & d
	\end{array} \right)=\left( \begin{array}{cc}
		a & -b \\
		-c & d
	\end{array} \right).$
	It is straight forward to see that $R$ is prime ring with involution $\ast$ and nontrivial endomorphism $g.$ Moreover, it is noticed that
	$Z(R)=\left\{ \left(\begin{array}{cc}
		a & 0 \\
		0 & a
	\end{array} \right) |~a\in R\right\}$
	and $X^{\ast}=X$ for all $X\in Z(R).$ It implies that $Z(R)\subseteq H(R),$ which implies that $\ast$ is the involution of the first kind. With $g_{1}=g$ and $g_{2}=g$, the assertions $(A_{1})$-$(A_{4})$ hold. However, $R$ is not commutative.
\end{exa}
In the following example illustrates that these results fail to be valid for semiprime rings.
\begin{exa}
	Let us consider $R=R_{1}\times R_{2},$ where $R_{1}=M_{2}(\mathbb{Q})$ and $R_{2}=\mathbb{C}.$ If $\ast$ is an involution defined on $R_{1}$ as in the previous example and $\dagger$ is an involution on $R$ defined as $(X,u)^{\dagger}=(X^{\ast},\overline{u}),$ which is an involution of the second kind on $R.$ If $g_{1},g_{2}$ are considered to be zero endomorphisms of $R,$ then it can be seen that all the conditions $(A_{1})$-$(A_{4})$ are satisfied, but $R$ is not commutative.
\end{exa}


\section*{Acknowledgements}
The authors express their gratitude to the learned referee and the editor for conducting an expert review which consequently enhanced the quality of this study.

\end{document}